\documentclass[10pt,draft]{article}

\usepackage[T1]{fontenc}
\usepackage{amsfonts}

\newtheorem{thm}{theorem}[section]
\newtheorem{theorem}[thm]{Theorem}
\newtheorem{proposition}[thm]{Proposition}
\newtheorem{lemma}[thm]{Lemma}
\newtheorem{corollary}[thm]{Corollary}
\newtheorem{remark}[thm]{Remark}

\newtheorem{definition}[thm]{Definition}

\begin{document}

\title{On the central polynomials with involution of $M_{1,1}(E)$}

\author{Diogo Diniz Pereira da Silva e Silva \footnote{Supported by CNPq}}
\date{{\small \textit{Unidade Acadêmica de Matemática e Estatística\\
Universidade Federal de Campina Grande\\
Cx. P. 10.044, 58429-970, Campina Grande, PB, Brazil}}\\
E-mail: diogo@dme.ufcg.edu.br}
\maketitle



\maketitle

\begin{abstract}
Let $K$ be an infinite field of characteristic $\neq 2$. In this paper we study the $*$-central polynomials for the algebra $M_{1,1}(E)$ with the involution $(*)$ induced by the transpose superinvolution on $M_{1,1}(K)$. More precisely, we determine a finite set of polynomials that together with the $*$-identities generate the $*$-central polynomials. In the case where the field $K$ has characteristic zero a finite set of generators for the $*$-central polynomials is determined. 

\end{abstract}

\section{Introduction}

The problem of the existence of central polynomials in a matrix algebra was raised by I. Kaplansky in \cite{Kaplansky} (see also \cite{Kaplansky2}). This problem was solved independently by Formanek \cite{Formanek} and Razmyslov \cite{Razmyslov}, and it has implications of interest in the study of polynomial identity rings (see \cite{Herstein}). Verbally prime algebras  were introduced by A. R. Kemer (see \cite{Kemer}) and play a proeminent role in the theory of PI-algebras. Razmyslov \cite{Razmyslov2} constructed  central polynomials for the algebras $M_{a,b}(E)$ and Kemer \cite{Kemer2} proved that every verbally prime algebra has a central polynomial. The existence of central polynomials for verbally prime algebras was used by A. Belov in \cite{Belov} to prove that each associative PI-algebra does not coincide with its commutant.  

Determining a basis for the space of all central polynomials of a given algebra is an important problem in the theory of PI-algebras. However generators for the central polynomials of an algebra are known in a few cases only. For the algebra $M_2(K)$ of $2\times 2$ matrices over a field $K$ generating sets were described in \cite{Okhitin} if $K$ is a field of characteristic zero and in \cite{ColomboKoshlukov2} for an infinite field of characteristic $p\neq 2$. In \cite{BrandaoPlamenAlexeiElida} the central polynomials of the Grassmann algebra of a vector space of countable dimension were described and over an infinite field of characteristic $p>2$ the central polynomials form a limit $T$-space. 

The study of polynomial identities naturally leads to some variants of the notion of polynomial identity (see \cite[Chapter 11]{BelovRowen}) such as weak identities, graded identities and identities with involution. This and the dificulty in determining basis in the ordinary case leads us to the study of these other types of identities and central polynomials. The verbally prime algebra $M_{1,1}(E)$ has been the subject of many research papers, see: \cite{DiVincenzoLaScala}, \cite{DiVincenzoLaScala2} for the weak identities; \cite{AzevedoeKoshlukov}, \cite{DiVincenzo} for the graded identities; \cite{DiVincenzoPlamen} for the identities with involution and \cite{BrandaoPlamen} for the graded central polynomials. 

The central polynomials with involution for $M_2(K)$ were described in \cite{BrandaoPlamen} and based on the results obtained in \cite{DiVincenzoPlamen} we study the central polynomials with involution for $M_{1,1}(E)$. The involution considered is the one induced by the transposition superinvolution on the superalgebra 
$M_2(K)$ with its natural $\mathbb{Z}_2$-grading. The $*$-identities of $M_n(E)$ and $M_{k,l}(E)$ with an involution induced by a pair of involutions on $M_n(K)$ and $E$ are also studied in \cite{DiVincenzoPlamen} and it is proved that, over a field of characteristic zero, the $*$-polynomial identities are just the ordinary ones. This implies that the 
$*$-central polynomials are also just the ordinary ones. In our main result we exhibit a finite set of polynomials that together with the identities generate the $*$-central polynomials of $M_{1,1}(E)$. 

\section{Central polynomials with involution}\label{sec2}

Let $K$ be an infinite field of characteristic $\neq 2$, in this article we consider algebras and vector spaces over $K$. Let $X$ be a countable set of indeterminates, the linear map $\#:K\langle X \rangle \rightarrow K\langle X \rangle$ determined by $(x_{i_1}\dots x_{i_k})^{\#}=(x_{i_k}\dots x_{i_1})$ is the \textit{reversal} involution on $K\langle X \rangle$. The \textit{canonical} involution $(*)$ on $K\langle X \rangle$ is the composition of the reversal involution with the endomorphism determined by $X_{2i-1}\mapsto X_{2i}$ and $X_{2i}\mapsto X_{2i-1}$ for all $i$. Henceforth $(K\langle X \rangle, *)$ denotes $K\langle X \rangle$ with the canonical involution, moreover to make the notation more convenient we write $x_i$ for $x_{2i-1}$ and $x_i^{*}$ for $x_{2i}$. The elements of this algebra will be called $*$-polynomials.

The algebra $(K\langle X \rangle, *)$ is the free algebra with involution, i. e., for any algebra with involution $(A,*)$ and for any $a_1,a_2,\dots$ in $A$ there exists a unique $*$-homomorphism $\varphi: (K\langle X \rangle, *) \rightarrow (A,*)$ such that $\varphi(x_i)= r_i$ for all $i$. An element $f(x_1,\dots,x_n,x_1^{*},\dots, x_n^{*})$ of 
$(K\langle X \rangle, *)$ is said to be a \textit{polynomial identity with involution} (or simply a $*$-identity) for $(A,*)$ if $f(a_1,\dots, a_n,a_1^{*},\dots, a_n^{*})=0$ for any $a_1,\dots, a_n $ in $A$. 

The set of $*$-identities of $(A,*)$ is denoted by $T_*(A)$, this is a two sided $*$-ideal $(K\langle X \rangle, *)$ wich is invariant under all its $*$-endomorphisms. Such $*$-ideals are called $T_*$-ideals. 

We denote by $C(A,*)$ the set \[\{f(x_1,\dots,x_n,x_1^{*},\dots, x_n^{*})\in (K\langle X \rangle, *)| f(a_1,\dots, a_n,a_1^{*},\dots, a_n^{*}) \in Z(A,*)\},\] this is a subspace of 
$(K\langle X \rangle, *)$ wich is invariant under all its $*$-endomorphisms. Such subspaces are called $*$-spaces. An element of $C(A,*)$ is called a \textit{central polynomial with involution} (of simply $*$-central polynomial) for $(A,*)$. The intersection of a family of $*$-subspaces of $(K\langle X \rangle, *)$ is also a $*$-subspace, hence given a set $S$ of polynomials we define $V(S)$ as the intersection of all $*$-subspaces that contain $S$ and we say that $V(S)$ is genereted as a $*$-space by $S$. 

Now we consider another set of generators for the free algebra with involution. Define\[y_i=\frac{x_i+x_i^{*}}{2}, z_i=\frac{x_i-x_i^{*}}{2}.\] We denote by $Y$ and $Z$ the sets $\{y_1,y_2,\dots,\}$ and $\{z_1,z_2,\dots,\}$. Note that the elements $y_i$ are symmetric and the elements $z_j$ are skewsymmetric.
Clearly $K\langle Y \cup Z\rangle = K\langle X \rangle$ and $f(y_1,\dots, y_m,z_1,\dots, z_n)$ is a $*$-central polynomial for $(A,*)$ (resp. $*$-identity) if and only if 
$f(a_1,\dots, a_m,b_1,\dots, b_n) \in Z(A,*)$ (resp. $=0$) for all $a_1,\dots, a_m \in A^{+}$ and $b_1,\dots, b_m \in A^{-}$.  

The degree of a polynomial $f$ in $K\langle Y \cup Z \rangle$ in the indeterminate $y_i$ (resp. $z_i$) is denoted by 
$deg_{y_i}(f)$ (resp. $deg_{z_i}(f)$). We say that $f$ is \textit{multihomogeneous} if every monomial of $f$ has the same degree in every indeterminate in $Y\cup Z$. Moreover if this degree is $1$ for every indeterminate that occurs in $f$ we say this polynomial is \textit{multilinear}.

Let $f$ be a polynomial in $ (K\langle X \rangle, *)$. Write it as $f=f_s+f_k$ where $f_s$ is the symmetric component of $f$ and $f_k$ is the skewsymmetric component.
\begin{remark}\cite[Remark 2.3.4]{RowenBook}
A $*$-polynomial $f(x_1,\dots,x_n,x_1^{*},\dots, x_n^{*})$ is a $*$-central polynomial but not a $*$-polynomial identity if and only if $[x_{n+1},f]$ and $f-f^{*}$, but not $f$, are 
$*$-polynomial identities.
\end{remark}
If $f$ is a $*$-central polynomial for $A$ the previous remark implies that $f_k=\frac{f^{*}-f}{2}$ is a $*$-polynomial identity for $A$. Thus $f_k \in C(A,*)$ and therefore $f_s=f-f_k$ lies in $C(A,*)$. This proves the following:

\begin{proposition}\label{P}
Let $f=f_s+f_k \in K\langle Y\cup Z \rangle$, where $f_s$ is symmetric and $f_k$ is skewsymmetric. If $f \in C(A,*)$ then $f_s \in C(A,*)$ and $f_k \in T_*(A)$.
\end{proposition}

Let $A$ be a ring, given $a_1,a_2$ we define the commutator of these two elements $[a_1,a_2]$ by $a_1a_2-a_2a_1$ and their Jordan product $a_1\circ a_2$ by $\frac{1}{2}(a_1a_2+a_2a_1)$. For $n>1$ and $a_1,\dots, a_n$ we define inductively $[a_1,\dots,a_n]=[[a_1,\dots, a_{n-1}],a_2]$. Denote by 
$B(Y\cup Z)$ the subalgebra of $K\langle Y\cup Z \rangle$ generated by $Z$ and the commutators $[a_1,\dots, a_n]$, where $n\geq 2$ and $a_i \in Y \cup Z$. The elements of $B(Y\cup Z)$ are called \textit{$*$-proper polynomials}.
\section{Involutions on $M_n(E)$ and $M_{k,l}(E)$ and $*$-Central Polynomials}\label{sec3}
In this section we present the results from \cite{DiVincenzoPlamen} that will be used in Sections \ref{sec4} and \ref{sec5}.

\subsection{Involutions on $M_n(E)$}

Let $A$ and $B$ be algebras with involutions $\alpha$ and $\beta$ respectively and $R$ the tensor product $A\otimes B$. The linear map $\gamma:R \rightarrow R$ determined by 
$\gamma(a\otimes b)=\alpha(a)\otimes \beta(b)$ is an involution on $R$. Recall that $E$ is the Grassmann algebra of the vector space $V$ with basis $\{e_1,e_2,\dots\}$. Let $\beta$  be a linear map of $V$ such that $\beta^{2}=1$. This map induces an involution on $E$, also denoted by $\beta$, such that $V$ is invariant under its action. Clearly every involution on $E$ such that $V$ is invariant is obtained in this way. Since we have the isomorphism $M_n(E) \cong M_n(K)\otimes E$ a pair $(\alpha, \beta)$ where $\alpha$ is an involution on $M_n(K)$ and $\beta$ is an involution on $E$ induces an involution $\gamma$ on $M_n(E)$.

\begin{proposition}\cite[Proposition 1]{DiVincenzoPlamen}
Let $K$ be a field of characteristic zero and $\gamma$ be the involution induced on $M_n(E)$ by a pair $(\alpha, \beta)$ of involutions defined on $M_n(K)$ and $E$,
respectively. If $V$ is invariant under the action of $\beta$ then any $*$-polynomial identity of $(M_n(E), \gamma )$ is trivial, that is
$f (y_1,\dots, y_l, z_1, \dots, z_m) \in T_{*}(M_n(E), \gamma )$ if and only if $f (x_1, \dots, x_{l+m}) \in T (M_n(E))$.
\end{proposition}

This immediatly yields the corresponding result for $*$-central polynomials.

\begin{corollary}\label{cccc}
Let $\gamma$ be the involution of the previous proposition. Then any $*$-central polynomial of $(M_n(E), \gamma )$ is trivial, that is
$f (y_1,\dots, y_l, z_1, \dots, z_m) \in C_{*}(M_n(E), \gamma )$ if and only if  $f (x_1, \dots, x_{l+m}) \in C (M_n(E))$.
\end{corollary}

\textit{Proof.}
If $f (y_1,\dots, y_l, z_1, \dots, z_m)$ is a $*$-central polynomial then the polynomial $[f (y_1,\dots, y_l, z_1, \dots, z_m),y_{l+1}]$ is an identity with involution. Hence the previous proposition implies that $[f (x_1,\dots, x_{l+m}),x_{l+m+1}]$ is an ordinary identity for $M_n(E)$ and therefore $f (x_1,\dots, x_{l+m})$ is an ordinary central polynomial.
\hfill $\Box$

The subalgebra $M_{k,l}(E)$ of $M_{k+l}(E)$ is the set of the matrices 
$\left(\begin{array}{cc}
	a&b\\
	c&d
\end{array}
 \right)$, where $a \in M_k(E_0)$, $b \in M_{k\times l}(E_1)$, $c \in M_{l\times k}(E_1)$, $d \in M_{l}(E_0)$ and $E_0$, $E_1$ are the subspaces of $E$ generated by the monomials of even length and the monomials odd  lenght respectively. Denote by $M_{k,l}(F)$ the superalgebra that consists of the algebra $M_{k+l}(F)$ with the elementary grading induced by the map 
$\mu_{k,l}:\{1,2,\dots, k+l\}\rightarrow \mathbb{Z}_2$ defined by $\mu_{k,l}(i)=0$ if $i\leq k$ and $\mu_{k,l}(i)=1$ otherwise. In this case the component $(M_{k+l}(F))_g$ is the subspace generated by the elementary matrices $e_{ij}$ such that $\mu(j)-\mu(i)=g$. The algebra $M_{k,l}(E)$ is isomorphic to $M_{k,l}(F)\widehat{\otimes} E $ and we consider involutions on $M_{k,l}(E)$ induced from involutions on $M_{k+l}(E)$. Let $\gamma$ be the involution on $M_n(E)$ by a pair $(\alpha, \beta)$ of involutions defined on $M_n(K)$ and $E$. Assume that:

\begin{enumerate}
\item[(1)] $\alpha$ preserves the $\mathbb{Z}_2$-graded components of $M_{k,l}(F)$;
\item[(2)] the vector space $V$ is invariant under the action of $\beta$.
\end{enumerate}

Under these hypothesis the involution $\gamma$ induces an involution on $M_{k,l}(E)$ and we have analogous results for the $*$-polynomial identities and $*$-central polynomials of 
$M_{k,l}(E)$.

\begin{proposition}\cite[Proposition 3]{DiVincenzoPlamen}
Let $K$ be a field of characteristic zero and $\gamma$ be the involution induced on $M_n(E)$ by a pair $(\alpha, \beta)$ of involutions defined on $M_n(K)$ and $E$,
respectively. If the involutions $\alpha$ and $\beta$ satisfy the conditions $(1)$ and $(2)$ above then $\gamma$ induces an ivolution on $M_{k,l}(E)$ and any $*$-polynomial identity of 
$(M_{k,l}(E), \gamma )$ is trivial, that is
$f (y_1,\dots, y_l, z_1, \dots, z_m) \in T_{*}(M_{k,l}(E), \gamma )$ if and only if  $f (x_1, \dots, x_{l+m}) \in T(M_{k,l}(E))$.
\end{proposition}

This immediatly yields the corresponding result for $*$-central polynomials. The proof is simple and analogous to the proof of Corollary \ref{cccc} and we omit it.

\begin{corollary}
Let $\gamma$ be the involution of the previous proposition. Then any $*$-central polynomial of 
$(M_{k,l}(E), \gamma )$ is trivial, that is
$f (y_1,\dots, y_l, z_1, \dots, z_m) \in C_{*}(M_{k,l}(E), \gamma )$ if and only if  $f (x_1, \dots, x_{l+m}) \in C(M_{k,l}(E))$.
\end{corollary}

\subsection{Superinvolutions and involutions on $M_{1,1}(E)$}

Let $A=A_0\oplus A_1$ be a superalgebra. A $\mathbb{Z}_2$-graded map $\circ: A \rightarrow A$ such that $(a^{\circ})^{\circ}=a$ for every $a$ in $A$ is a superinvolution if 
\[(ab)^{\circ}=(-1)^{||a||\cdot ||b||}b^{\circ}a^{\circ}\] for every homogeneous elements $a$, $b$ of degree $||a||$, $||b||$ respectively in the $\mathbb{Z}_2$-grading of $A$. Superinvolutions play an important role in ring theory (see \cite{GomezShestakov}). If $A=A_0\oplus A_1$ and $B=B_0\oplus B_1$ are superalgebras with superinvolutions $\circ$ and $\diamond$ then $R=A\widehat{\otimes}B = A_0\otimes B_0 \oplus A_1 \otimes B_1$ has an involution $*$ defined by $(a\otimes b)^{*}=a^{\circ}\otimes b^{\diamond}$. Since the identity map is a superinvolution on $E$ and $M_{k,l}(E)\cong M_{k,l}\widehat{\otimes} E$ any superinvolution on the superalgebra $M_{k,l}(F)$ induces an involution on $M_{k,l}(E)$. For $k=l=1$ we have the following result:

\begin{theorem}\cite[Theorem 3.2]{GomezShestakov}
The only two superinvolutions on $M_{1,1}(F)$ are $trp$ and $(trp)p$,
where $p$ is the automorphism of $M_{1,1}(F)$ given by $p(a_0+a_1)=a_0-a_1$ (the
parity automorphism) and \[trp\left(\begin{array}{cc}
	a & b \\
	c & d
\end{array}\right)=\left(\begin{array}{cc}
	d & b \\
	-c & a
\end{array}\right).\]
\end{theorem}

Let $*$ and $\odot$ be the involutions on $M_{1,1}(E)$ induced by the superinvolutions $trp$ and $(trp)p$ respectively The map $\phi\left( \left(\begin{array}{cc}
	a & b \\
	c & d
\end{array}\right) \right)= \left(\begin{array}{cc}
	d & c \\
	b & a
\end{array}\right)$ is an isomorphism of the algebras with involution $((M_{1,1})(E), *)$ and $(M_{1,1}(E), \odot)$ (see \cite{DiVincenzoPlamen}). Hence the $*$-central polynomials (and the $*$-polynomial identites) of these two algebras with involution are the same.

We denote by $(R,*)$ (or simply by $R$) the algebra \[M_{1,1}(E)=\left\{ 
\left(\begin{array}{cc}
	a & b \\
	c & d
\end{array}\right)|a, d \in E_0\mbox{, } b,c \in E_1
 \right\},\] with the involution $*$ defined by \[\left(\begin{array}{cc}
	a & b \\
	c & d
\end{array}\right)^{*}=\left(\begin{array}{cc}
	d & b \\
	-c & a
\end{array}\right).\] 

The subspaces of symmetric and skewsymmetric elements are
 \[R^{+}=\left\{ 
\left(\begin{array}{cc}
	a & b \\
	0 & a
\end{array}\right)|a \in E_0\mbox{, } b \in E_1
 \right\}\] and

\[R^{-}=\left\{ 
\left(\begin{array}{cc}
	a & 0 \\
	b & -a
\end{array}\right)|a \in E_0\mbox{, } b \in E_1
 \right\},\] respectively.

In the next proposition we list some $*$-polynomial identities for $R$.
\begin{proposition}\cite[Proposition 4]{DiVincenzoPlamen} \label{Id}
The following polynomials lie in the ideal $T_*(R)$:
\begin{eqnarray*}
 H_1&=&[y_1,y_2],\\
 H_2&=&z_1z_2z_3-z_3z_2z_1,\\
 H_3&=&[z_1,z_2][z_3,z_4],\\
 H_4&=&[y_1,z_1,z_2,z_3]-2(z_1\circ z_2)[y_1,z_3]),\\
 H_5&=&[z_1,z_2]y_1[z_3,z_4]+[z_3,z_4]y_1[z_1,z_2],\\
 H_6&=&2z_1[y_1,z_2,z_3]+[y_1,z_1,z_2,z_3]+[z_1,z_2][y_1,z_3]+\\&&[z_1,z_3][y_1,z_2]+[z_2,z_3][y_1,z_1],\\
 H_7&=&2[z_1,z_2]z_3[y_1,z_4]+[z_1,z_2][y_1,z_3,z_4],\\
 H_8&=&[y_1,z_1][y_2,z_2]+[y_1,z_2][y_2,z_1],\\
 H_9&=&[y_1,z_1][y_2,z_2,z_3]-[y_1,z_2,z_1][y_2,z_3],\\
 H_{10}&=&[z_1,z_2][y_1,z_3][y_2,z_4]-[z_1,z_4][y_1,z_2][y_2,z_3]-\\&&[z_2,z_3][y_1,z_1][y_2,z_4]+[z_3,z_4][y_1,z_1][y_2,z_2].
\end{eqnarray*}
\end{proposition}

Note that the above proposition holds for an infinite field $K$. As a consequence we conclude that the following polynomials lie in $C(R,*)$:

\begin{enumerate}
\item[(a)] $z_1\circ z_2$;
\item[(b)] $[y_1,z_1,z_2]-2(z_1\circ z_2)y_1$;
\item[(c)] $y_1^{p}$, if $\textit{char}K=p> 2$.
\end{enumerate}

We denote by $I$ the $*$-ideal $T_*(R)$ and by $V$ the $*$-space in $K\langle Y \cup Z \rangle$ generated by $I$ and by the polynomials $(a)$ and $(b)$ above together with the polynomial $(c)$ if $\textit{char}K=p> 2$.

\section{$*$-Proper Central Polynomials of $M_{1,1}(E)$}\label{sec4}

In this section we study the $*$-proper central polynomials of $R$. In this section and in Section \ref{sec5} the field $K$ is infinite.

\begin{proposition}\label{Pp}
If $f(y_1,\dots, y_l,z_1,\dots,z_m)$ is a $*$-proper central polynomial for $R$ that depends on the variables in $Y$ then $f$ is equivalent modulo $I$ to a linear combination of the polynomials

\[[y_{a_1},z_{i_1},\dots,z_{i_k},z_{j_1}][y_{a_2},z_{j_2}]\dots[y_{a_l},z_{j_l}],\] 
where $a_1<\dots < a_l$, $i_1\leq i_2\leq \dots \leq i_k$, $j_1<\dots<j_l$.

\end{proposition}

\textit{Proof.}
As a consequence of Proposition \ref{P} we may assume that $f$ is symmetrical. Hence the result follows from \cite[Lemma 13]{DiVincenzoPlamen}, 
\cite[Proposition 15]{DiVincenzoPlamen} and \cite[Proposition 16]{DiVincenzoPlamen}.
\hfill $\Box$

Next we obtain some consequences of the generators of $V$, the following remark will be usefull.

\begin{remark}\cite[Remark 6, Lemma 7]{DiVincenzoPlamen}\label{R}
Let $l\geq 1$. For any $\sigma, \tau \in S_l$ the polynomials
\[ [y_1,z_1]\dots [y_l,z_l]-(-1)^{\sigma \tau}[y_{\sigma(1)},z_{\tau(1)}]\dots [y_{\sigma(l)},z_{\tau(1)}],\] and 
\[[y_1,z_{\sigma(1)}, \dots, z_{\sigma(l)}, z_{l+1}]-[y_1,z_1,\dots, z_l, z_{l+1}]\] lie in $I$.
\end{remark}

Let $l\geq 1$ be a natural number, we consider the polynomials \[P_1(y_2,\dots, y_{l+1},z_1,\dots, z_{l+1})=[y_2,z_1,z_2][y_3,z_3]\dots [y_{l+1},z_{l+1}]\] and 
\[P_i(y_2,\dots, y_{l+1},z_1,\dots, z_{l+1})=[y_2,z_i,z_1]\dots [y_{i},z_{i-1}][y_{i+1},z_{i+1}]\dots [y_{l+1},z_{l+1}],\] for $1<i\leq l+1$. 

By induction on $l$ it follows from the previous remark and the identity $H_9$ in Proposition \ref{Id} that 

\begin{equation}\label{eee}
[P_i,z_{l+2}]=\sum_{j<i}(-1)^{l-j}P_{ji}+\sum_{j>i} (-1)^{l+1-j}P_{ij},
\end{equation}
where \[P_{m,n}=[y_{2},z_{m},z_{n},z_{k_1}][y_{3},z_{k_2}]\dots[y_{l+1},z_{k_{l}}],\] the indexes are ordered $k_1<\dots<k_{l}$, $i<j$ and
 $\{m,n\} \cup\{k_1,\dots, k_{l}\}$ is a partition of $\{1,2,\dots,l+2\}$. Note that (\ref{eee}) also holds for $i=1$.

For $l\geq 1$ denote
\begin{equation}\label{Ed}
C_l=\sum_{n=1}^{l+1} (-1)^{n} P_n,\mbox{ and } D_l=y_1C_l-[y_1,z_1]\dots[y_{l+1},z_{l+1}].
\end{equation}
It follows from (\ref{eee}) that $[C_l,z_{l+2}]$ and $[D_l,z_{l+2}]$ lie in $I$. Since $[y_1,y_2]$ it follows that $[C_l,y_{l+2}]$ and $[D_l,y_{l+2}]$ are $*$-identities. Hence we conclude that $C_l$ and $D_l$ are $*$-central polynomials for $R$. In the next lemma we prove that in fact $C_l,D_l \in V$.

\begin{lemma}\label{L}
The $*$-space $V$ contains the polynomials \[G_n=(z_1\circ z_2)(z_3\circ z_4)\dots (z_{2n-1}\circ z_{2n}),\] and the polynomials $C_l$,$D_l$ defined in (\ref{Ed}).
\end{lemma}

\textit{Proof.}
Clearly $G_1$ lies in $V$ and it follows by induction on $n$ that $G_n$ lies in $V$ for every $n$. It remains to prove that $C_l$ and $D_l$ lie in $V$. The polynomial 

\begin{equation}\label{Ec}
C_1=[y_2,z_2,z_1]-[y_2,z_1,z_2],
\end{equation}
is a consequence of the polynomial $(b)$. Note that we may substitute $z_1$ by $(y_1\circ z_1)$ in (\ref{Ec}), thus we conclude that

\begin{equation}\label{E1}
D_1=y_1([y_2,z_2,z_1]-[y_2,z_1,z_2])+[y_1,z_1][y_2,z_2],
\end{equation}
lies in V.
For $l>1$ we substitute $y_1$ by $Y=[y_3,z_3]\dots [y_{l+1},z_{l+1}]$ in (\ref{E1}), thus

\begin{equation}\label{Cl}
Y([y_2,z_2,z_1]-[y_2,z_1,z_2])+[Y,z_1][y_2,z_2] \in V.
\end{equation}

Note that two symmetric polynomials commute modulo $I$, therefore we have

\[Y([y_2,z_2,z_1]-[y_2,z_1,z_2])\equiv_I P_2-P_1.\] Moreover $[Y,z_1]=\sum_{i=3}^{l+1}[y_{3},z_3]\dots [y_{i},z_i,z_1]\dots [y_{l+1},z_{l+1}]$ and it follows from Remark \ref{R} that 

\[[Y,z_1][y_2,z_2]\equiv_I \sum_{n=3}^{l} (-1)^{n}P_n,\] therefore the polynomial in (\ref{Cl}) is $C_l$.  Substituting $z_1$ by $(y_2\circ z_1)$  in $C_l$ we conclude that $D_l$ lies in $V$.

\hfill $\Box$

We denote by $W$ the subspace of $V$ generated as a $*$-space by the polynomials $G_n$, $C_l$ and $C_lG_n$, where $n>0$ and $l>0$. Our main objective now is to prove Corollary \ref{C}, the proof is divided in the next three lemmas.

\begin{lemma}\label{Lz}
If $f(z_1,\dots,z_m)$ is a non-constant central polynomial for $R$ then $f$ is equivalent modulo $I$ to a polynomial in the $*$-space generated by the polynomials $G_n$. In particular 
$f$ is equivalent modulo $I$ to a polynomial in $W$.
\end{lemma}

\textit{Proof.}
Since the field is infinite we may assume that $f$ is multihomogeneous. Note that $z_1z_2z_3-z_3z_2z_1$ is a $*$-polynomial identity for $R$, hence $z_{i_1}\dots z_{i_{2n+1}}$ is equivalent modulo $I$ to a skewsymmetric polynomial. Therefore if $f$ has odd degree we have $f(z_1,\dots,z_m)\equiv_I f_1(z_1,\dots,z_m)$ where $f_1$ is a skewsymmetric polynomial. By Proposition \ref{P} we conclude that $f_1$ lies in $I$. Now assume that $f$ has degree $2n$. Using the equality 
$z_{i_1}z_{i_2}=z_{i_1}\circ z_{i_2}+\frac{1}{2}[z_{i_1},z_{i_2}]$ and the $*$-polynomial identity $[z_1,z_2][z_3,z_4]$ we get 
\begin{equation}\label{e22}
z_{i_1}z_{i_2}\dots z_{i_{2n-1}}z_{i_{2n}}\equiv_{I} (z_{i_1}\circ z_{i_2})\dots (z_{i_{2n-1}}\circ z_{i_{2n}})+g,
\end{equation}
where 
\[g=\sum_{k=1}^{n}\frac{1}{2}(z_{i_1}\circ z_{i_2})\dots \widehat{(z_{i_{2k-1}}\circ z_{i_{2k}})}\dots (z_{i_{2n-1}}\circ z_{i_{2n}})[z_{i_{2k-1}},z_{i_{2k}}],\] and the symbol 
\hspace{0.3cm}$\widehat{ }$\hspace{0.3cm} over the polynomial $(z_{i_{2k-1}}\circ z_{i_{2k}})$ means it does not occur in the product. Since $g^{*}\equiv_I -g$ we have 
$g\equiv_I \frac{1}{2}(g-g^{*})$ and we may assume that the polynomial $g$ in (\ref{e22}) is skewsymmetric. Hence $f\equiv_I f_1 + g_1$, where $f_1$ is in the $*$-space generated by 
$G_n$ and $g_1$ is a skewsymmetric polynomial. Since $f$ is central it follows that $g_1$ is central and Proposition \ref{P} implies that $g_1 \in I$.
\hfill $\Box$

The following lemma deals with $*$-central polynomials in $R$ that depend on the variables in $Y$ only. 

\begin{lemma}\label{L4}
Let $g(y_1,\dots,y_m)$ be a non-constant central polynomial for $R$. We have:

\begin{enumerate}
\item[(i)] $g \in I$ if $char K =0$;
\item[(ii)] $g\equiv_I \alpha (y_1^{a_1p})\dots (y_m^{a_mp}) \in V$ if $char K =p>2$.
\end{enumerate}
\end{lemma}

\textit{Proof.}
We consider 
$\overline{y_i}=\left(\begin{array}{cc}
\alpha_i & \beta_i\\
0 & \alpha_i	
\end{array}\right),$ where $\alpha_i \in E_0$ and $\beta_i \in E_1$. We have $g(\overline{y_1},\dots, \overline{y_m})=\lambda \overline{y_1}^{\alpha_1}\dots \overline{y_m}^{\alpha_m}$. Since \[\overline{y_i}^{k_i}=\left(\begin{array}{cc}
\alpha_i^{k_i} & k_i\alpha_i^{k_i-1}\beta_i\\
0 & \alpha_i^{k_i}	
\end{array}\right),\] the $(1,2)$-entry of $\overline{y_1}^{\alpha_1}\dots \overline{y_m}^{\alpha_m}$ is 
\[m=\sum_{i=1}^{m}k_i\beta_i\alpha_1^{k_1}\dots \alpha_{i-1}^{k_{i-1}}\alpha_{i}^{k_{i}-1}\alpha_{i+1}^{k_{i+1}}\dots \alpha_{m}^{k_{m}}.\] Thus if $char K = 0$ it follows that 
$\lambda =0$ and we have $(i)$ and if $char K =p >2$ we conclude that if $\lambda \neq 0$ each $k_i$ is divisible by $p$ and we have $(ii)$.
\hfill $\Box$

\begin{lemma}\label{L1}
If $f(y_1,\dots, y_l,z_1,\dots,z_m)$ is a $*$-proper central polynomial for $R$ that depends on the variables in $Y$ and in $Z$ then $f$ is equivalent modulo $I$ to a polynomial in the $*$-space $W_0$ generated by the polynomials $C_l$. In particular $f$ is equivalent modulo $I$ to a polynomial in $W$.
\end{lemma}

\textit{Proof.}
Since $K$ is an infinite field we may assume that $f$ is multihomogeneous. Denote by $\textbf{J}$ the set of $l$-tuples $\textbf{j}=(j_1,\dots, j_l)$ with $1\leq j_1<\dots< j_l\leq m$ such that $z_{j_1}, \dots, z_{j_l}$ appear in $f$. Let $\textbf{j}=(j_1,\dots, j_l)$ we denote by $P_{\textbf{j}}$ the polynomial 
\[[y_{1},z_{i_1},\dots,z_{i_k},z_{j_1}][y_{l},z_{j_2}]\dots[y_{l},z_{j_l}],\] of the same multidegree as $f$. It follows from Proposition \ref{Pp} that 
\[f\equiv_I\sum_{\textbf{j}\in \textbf{J}}\alpha_jP_{\textbf{j}}.\] We may assume that $z_m$ appears in $f$ and since $C_l$ lies in $W_0$ we conclude that each $P_{\textbf{j}}$ is equivalent modulo $W_0$ to a linear combination of polynomials \[P_{\textbf{j}^{\prime}}=[y_{1},z_{i_1},\dots,z_{i_k},z_{j_1}][y_{l},z_{j_2}]\dots[y_{l},z_{m}],\] with the commutator 
$[y_l,z_m]$ as the last factor in the product. Hence there exists a $*$-proper polynomial $g$ such that $f-g[y_l,z_m]$ lies in $W_0$. In this case $g[y_l,z_m]$ is a central polynomial for 
$R$ and this implies that $g$ lies in $I$.
\hfill $\Box$

\begin{corollary}\label{C}
If $w(y_1,\dots,y_l,z_1,\dots,z_m)$ is a $*$-proper central polynomial for $R$ and $w$ is not a constant then $w$ is equivalent modulo $I$ to a polynomial in $W$. In particular we conclude that $w$ lies in $V$.
\end{corollary}

\textit{Proof.}
The result follows from lemmas \ref{Lz}, \ref{L1} and \ref{L4}.
\hfill $\Box$

\section{$*$-Central Polynomials of $M_{1,1}(E)$}\label{sec5}

It follows from the Poicaré-Birkhoff-Witt theorem that every polynomial \\ $f(y_1,\dots, y_l,z_1,\dots,z_m)$ in $K\langle Y \cup Z \rangle$ may be uniquely written as 

\begin{equation}\label{E}
f=\sum y_1^{a_1}\dots y_{l}^{a_l} w_{\textbf{a}},
\end{equation}
where
$\textbf{a}=(a_1,\dots,a_n)$ is an $n$-tuple of natural numbers and $w_{\textbf{a}}$ is a $*$-proper polynomial. 

We recall the definition of the rank of a multihomogeneous polynomial in \cite{BrandaoPlamen}. 

\begin{definition}
Let $f(y_1,\dots, y_l,z_1,\dots,z_m) \neq 0$ be as in (\ref{E}). We define the \textbf{rank} of $f$, denoted by $r(f)$, as the greatest $n$-tuple $\textbf{a}=(a_1,\dots, a_n)$, in the lexicographical order, such that $w_{\textbf{a}}\neq 0$.
\end{definition}

\begin{remark}
Let $A$ be an algebra with involution. It is known that if we write $f$ as in ($\ref{E}$) then $f$ lies in $T_*(A)$ if and only if each $w_{\textbf{a}}$ lies in $T_*(A)$ (see 
\cite[Lemma 2.1]{DrenskyGiambruno}). This is not the case for $*$-central polynomials, a simple counterexample is the central polynomial $(b)$ of $R$. However we have the following proposition.
\end{remark}

\begin{proposition}\label{Pw}
Let $f(y_1,\dots, y_l,z_1,\dots,z_m)$ be as in (\ref{E}). If $f$ is a $*$-central polynomial of rank $\textbf{a}$ for an algebra with involution $A$ then the polynomial $w_{\textbf{a}}$ is a $*$-central polynomial for $A$.
\end{proposition}
\textit{Proof.}
The proof of this result is a part of the proof of \cite[Theorem 21]{BrandaoPlamen}.
\hfill $\Box$

The main part in the proof of Theorem \ref{T} is to prove that $C(R,*)\subset V$. The next result allows us to use the notion of rank of a polynomial to reduce the proof of this inclusion to the case of $*$-proper polynomials.

\begin{lemma}\label{L3}
Let $w(y_{1},\dots, y_l,z_1,\dots, z_m) \in W$ be a $*$-proper polynomial and  $y_{1}^{a_1}\dots y_{l}^{a_l}$ be a monomial. If the product $y_{1}^{a_1}\dots y_{l}^{a_l} w$ does not lie in $V$ then it is equivalent modulo $V$ to a polynomial of rank $< \textbf{a}=(a_1,\dots, a_l)$.
\end{lemma}
\textit{Proof.}
If $w$ depends on the variables $z_i$ only then it follows from Lemma \ref{Lz} that we only need to consider the case $w=G_n$. Note that $(z_1\circ z_2)y_1-\frac{1}{2}[y_1,z_1,z_2]$ lies in $V$, and by induction on $n$ it follows that 
\begin{equation}\label{E0}
G_ny_1-\frac{1}{2^n}[y_1,z_1,z_2,\dots, z_{2n}],
\end{equation}
also lies in $V$. We substitute $y_1$ by $y_{1}^{a_1}\dots y_{l}^{a_l}$ in this polynomial and conclude that $y_{1}^{a_1}\dots y_{l}^{a_l}w$ is equivalent modulo $V$ to 
\[f=\frac{1}{2^n}[y_{1}^{a_1}\dots y_{l}^{a_l},z_1,z_2,\dots, z_{2n}].\] We have 
\begin{equation}\label{l22}
[y_{1}^{a_1}\dots y_{l}^{a_l},z_1]\equiv_I \sum_i a_iy_1^{a_1}\dots y_{i-1}^{a_{i-1}} y_i^{a_i-1} y_{i+1}^{a_{i+1}} \dots y_l^{a_l}[y_i,z_j].
\end{equation}
Since $y_{1}^{a_1}\dots y_{l}^{a_l} w$ does not lie in $V$ it follows that $f$ is not a $*$-identity. Hence using (\ref{l22}) it follows from induction on $n$ that it is equivalent modulo $I$ to a polynomial of rank $<\textbf{a}$. Suppose now that $w$ depends on the variables in $Y$ and in $Z$. In this case it follows from Lemma \ref{L1} that we only need to consider the case $w=C_l$. Substitute $y_1$ by 
$y_{1}^{a_1}\dots y_{l}^{a_l}$ in 
\begin{equation}\label{EE1}
D_l=y_1C_l-[y_1,z_1]\dots[y_{l+1},z_{l+1}].
\end{equation}
Since $D_l$ lies in $V$ we get that $y_{1}^{a_1}\dots y_{l}^{a_l}C_l$ is equivalent modulo $V$ to 
\[f=[y_{1}^{a_1}\dots y_{l}^{a_l},z_1]\dots[y_{l+1},z_{l+1}].\] As in the previous case $f$ is equivalent modulo $I$ to a polynomial of rank $<\textbf{a}$. 

\hfill $\Box$

\begin{theorem}\label{T}
The $*$-space $C(R,*)$ is generated by $T_*(R)$ together with the polynomials $z_1\circ z_2$, $[y_1,z_1,z_2]-[y_1,z_2,z_1]$ and $y_1^{p}$ if $\textit{char}K=p$. In other words 
$C(R,*)=V$.
\end{theorem}

\textit{Proof.}
We already have the inclusion $V\subset C(R,*)$, therefore we need only to prove the reverse inclusion, since the field is infinite we consider multihomogeneous polynomials only. 
Let $f$ be a $*$-central polynomial for $R$ of rank $\textbf{a}=(a_1,\dots, a_n)$ and write $f$ as in (\ref{E}). It follows from Proposition \ref{Pw} that $w_{\textbf{a}}$ is a $*$-central polynomial for $R$. In this case Corollary \ref{C} implies that $w_a$ lies in $W$ and using Lemma  \ref{L3} we conclude that either 
$y_1^{a_1}\dots y_{l}^{a_l} w_{\textbf{a}} \in V$ or $y_1^{a_1}\dots y_{l}^{a_l} w_{\textbf{a}}$ (and hence $f$) is equivalent modulo $V$ to a polynomial of rank $< \textbf{a}$. Repeating this argument if necessary we conclude that $f$ is equivalent modulo $V$ to a $*$-proper central polynomial, hence the result follows from Corollary \ref{C}.
\hfill $\Box$

As a consequence of the previous theorem the problem o determining a basis for the central polynomials of $R$ is a consequence of determining a basis, as a $*$-ideal, for the identities of $R$. In particular we have the following corollary.

\begin{corollary}
 If $\textit{char}K=0$ the $*$-space $C(R,*)$ is generated by the polynomials $z_1\circ z_2$, $[y_1,z_1,z_2]-[y_1,z_2,z_1]$ and the polynomials $(x_1H_1x_2), \dots, (x_1H_{10}x_2)$.
\end{corollary}
\textit{Proof.}
The result follows from Theorem \ref{T} and \cite[Theorem 5]{DiVincenzoPlamen}.
\hfill $\Box$

\begin{flushleft}
\textbf{Acknowledgements}
\end{flushleft}
This work was supported by CNPq under Grant 303534/2013-3.

\end{document}